# A New Formulation for Total Least Square Error Method in *d*-dimensional Space with Mapping to a Parametric Line

Vaclav Skala

*Department of Computer Science and Engineering, Faculty of Applied Sciences, University of West Bohemia, Univerzitni 8, CZ 306 14 Plzen, Czech Republic*
*http://www.VaclavSkala.eu*

**Abstract.** There are many practical applications based on the Least Square Error (LSE) or Total Least Square Error (TLSE) methods. Usually the standard least square error is used due to its simplicity, but it is not an optimal solution, as it does not optimize distance, but square of a distance. The TLSE method, respecting the orthogonality of a distance measurement, is computed in $d$-dimensional space, i.e. for points given in $E^2$ a line $\pi$ in $E^2$, resp. for points given in $E^3$ a plane $\rho$ in $E^3$, fitting the TLSE criteria are found. However, some tasks in physical sciences lead to a slightly different problem.

In this paper, a new TSLE method is introduced for solving a problem when data are given in $E^3$ a line $\pi \in E^3$ is to be found fitting the TLSE criterion. The presented approach is applicable for a general $d$-dimensional case, i.e. when points are given in $E^d$ a line $\pi \in E^d$ is to be found. This formulation is *different* from the TLSE formulation.

**Keywords:** Least square error, total least square error, approximation, optimality
**PACS:** 02.60.-x, 02.30.Jr, 02.60 Dc

## INTRODUCTION

In many applications interpolation and approximation methods are used. Applications can be found in solution of economical, statistical and technical problems. Sometimes such techniques are called as “Regression” methods with several attributes, like linear, quadratic etc. The “standard” Least Square Error (LSE) can be found in many publications and it is based actually on differences evaluation on the dependent axis, i.e. on $y$-axis if $y = f(x)$ is considered. The solution leads to pseudo-inverse of a rectangular matrix as we have over determined system of linear equations [1][2][5][12]. However, if the problem is formulated implicitly as $F(x, y) = 0$ this approach is not directly applicable, see Fig.1 for a difference.

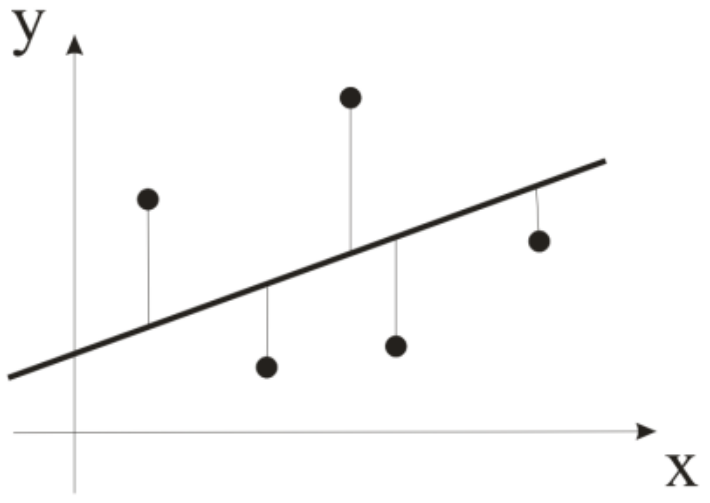

Fig. 1.a: Standard LSE for $y = f(x)$

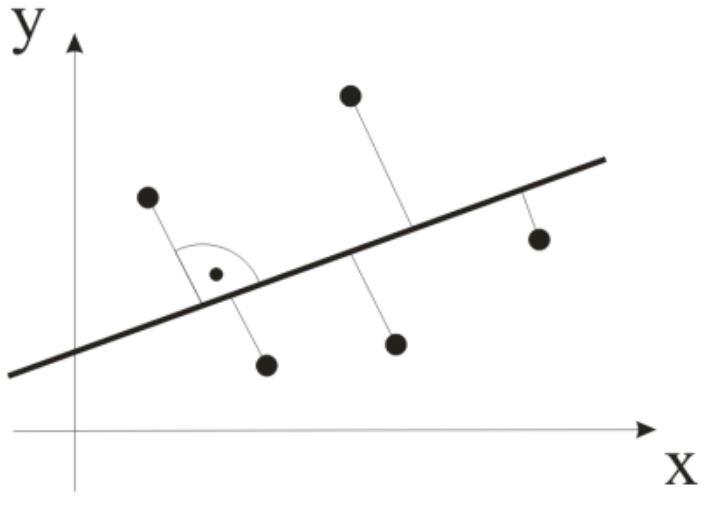

Fig. 1.b: Total LSE for $F(x, y) = 0$

It can be seen, that there is a significant difference between those two cases. However, especially in geometrical problems, the TSLE computation is not taken and received results can be far from the optimal one.

## LEAST SQUARE ERROR APPROXIMATION

Let us consider the “standard” LSE approximation for an explicit function [5][6]:

$$y = \boldsymbol{\xi}^T \boldsymbol{x} + \delta \tag{1}$$

where $\boldsymbol{x}$ is n-dimensional value, $\boldsymbol{\xi}$ is a vector and $\delta$ are unknown coefficients in general, e.g. of a line or plane etc. As we have more given points than unknown coefficients we get an over determined system of linear equations:

$$\boldsymbol{A\xi} = \boldsymbol{b} \tag{2}$$

where the matrix $\boldsymbol{A}$ contains the given points $\boldsymbol{x}$ and $\boldsymbol{b}$ contains $y$ value for the each given point.

The error is then defined as

$$r = \|\boldsymbol{Ax} - \boldsymbol{b}\| \tag{3}$$

then

$$\begin{aligned} r^2 &= (\boldsymbol{Ax} - \boldsymbol{b})^T(\boldsymbol{Ax} - \boldsymbol{b}) \\ &= (\boldsymbol{Ax})^T(\boldsymbol{Ax}) - (\boldsymbol{Ax})^T\boldsymbol{b} - \boldsymbol{b}^T\boldsymbol{Ax} + \boldsymbol{b}^T\boldsymbol{b} \end{aligned} \tag{4}$$

To find an extreme, the following condition must be valid

$$\frac{\partial r^2}{\partial \boldsymbol{x}} = 2\boldsymbol{A}^T\boldsymbol{Ax} - \boldsymbol{A}^T\boldsymbol{b} - \boldsymbol{b}^T\boldsymbol{A} = 0 \tag{5}$$

which leads to a system of linear equations if $\boldsymbol{A}$ is symmetric.

$$\boldsymbol{A}^T\boldsymbol{Ax} = \boldsymbol{b}^T\boldsymbol{A} \tag{6}$$

The matrix $\boldsymbol{A}^T\boldsymbol{A}$ is a symmetric positively semi-definite matrix.

It should be noted that direct solution using pseudo-inverse is not convenient for high dimensional case as the system of linear equations is ill conditioned [9][10], in general.

In the following we explore the case, when Total Least Square Error (TLSE) method is requited but for a case, when points are given in $E^3$ and a line $\pi \in E^3$, fitting the TLSE criterion, is to be found. Note that this problem formulation is different from the TLSE formulation.

## TOTAL LEAST SQUARE ERROR APPROXIMATION

Let us consider a problem, when the given points $\left\{\langle \boldsymbol{x}_p \rangle\right\}_{p=1}^{N}$, $\boldsymbol{x}_p \in E^d$ and $d$-dimensional Euclidean space.

We want to find a line $\pi$ in $d$-dimensional space having a minimal orthogonal distance from all the given points.

For simplicity, we set $d = 3$ without any limitation to the $d$-dimensional case, see Fig.2. If the set of given points is extended, for simplicity of derivation (is not actually needed), as:

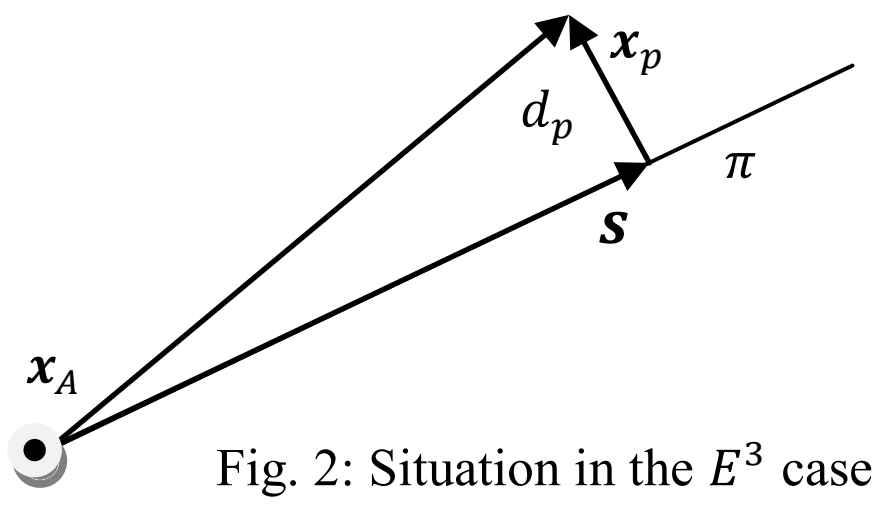

Fig. 2: Situation in the $E^3$ case

$$\boldsymbol{\Phi} = \left\{\langle \boldsymbol{x}_p \rangle\right\}_{p=1}^{N} \cup \left\{\langle -\boldsymbol{x}_p \rangle\right\}_{p=1}^{N} \tag{7}$$

then a line $\pi$, we are looking for, has to pass the origin of the coordinate system as:

$$\boldsymbol{x}_T = \frac{1}{2N}\sum_{p=1}^{N}\left(\boldsymbol{x}_p + \boldsymbol{x}_{N+p}\right) = \frac{1}{2N}\sum_{p=1}^{N}\left(\boldsymbol{x}_p - \boldsymbol{x}_p\right) = \boldsymbol{0} \tag{8}$$

as $\boldsymbol{x}_p = -\boldsymbol{x}_{N+p}$ ; the original data set of given points was extended. The given set of points $\boldsymbol{\Phi}$ can be used directly, if the given points are shifted, i.e. $\boldsymbol{x}_T = \boldsymbol{0}$.

As area $A$ of the triangle, see Fig.2, is given as:

$$2A = \left\|\boldsymbol{s} \times \left(\boldsymbol{x}_p - \boldsymbol{x}_A\right)\right\| = \|\boldsymbol{s}\|\, d_p \tag{9}$$

then the square of the distance $d_p^2$ is given as:

$$d_p^2 = \frac{\left(\boldsymbol{s} \times \left(\boldsymbol{x}_p - \boldsymbol{x}_A\right)\right)^T \left(\boldsymbol{s} \times \left(\boldsymbol{x}_p - \boldsymbol{x}_A\right)\right)}{\boldsymbol{s}^T\boldsymbol{s}} \tag{10}$$

As $\boldsymbol{x}_T = \boldsymbol{x}_A = \boldsymbol{0}$, we get a total orthogonal distance $D$ for all the given points as:

$$D = \sum_{p=1}^{2N} d_p^2 = \sum_{p=1}^{2N} = \frac{\left(\boldsymbol{s} \times \boldsymbol{x}_p\right)^T\left(\boldsymbol{s} \times \boldsymbol{x}_p\right)}{\boldsymbol{s}^T\boldsymbol{s}} = \sum_{p=1}^{2N} \frac{\left(\boldsymbol{x}_p \times \boldsymbol{s}\right)^T\left(\boldsymbol{x}_p \times \boldsymbol{s}\right)}{\boldsymbol{s}^T\boldsymbol{s}} \tag{11}$$

Now, we need to use matrix-vector operations, instead of the cross product.

A cross product $\boldsymbol{a} \times \boldsymbol{s}$ can be rewritten as:

$$\boldsymbol{a} \times \boldsymbol{s} = \begin{bmatrix} 0 & -a_z & a_y \\ a_z & 0 & -a_x \\ -a_y & a_x & 0 \end{bmatrix} \boldsymbol{s} = \boldsymbol{Ts} \tag{12}$$

We can write using matrix form of the cross product using a matrix $\boldsymbol{Q}_p$:

$$\boldsymbol{x_p} \times \boldsymbol{s} = \boldsymbol{Q}_p\boldsymbol{s} = \begin{bmatrix} 0 & -z_p & y_p \\ z_p & 0 & -x_p \\ -y_p & x_p & 0 \end{bmatrix} \boldsymbol{s} \tag{13}$$

Multiplication of a transpose matrix $\boldsymbol{T}^T$ and matrix $\boldsymbol{T}$ leads to:

$$\boldsymbol{T}^T\boldsymbol{T} = \begin{bmatrix} 0 & a_z & -a_y \\ -a_z & 0 & a_x \\ a_y & -a_x & 0 \end{bmatrix} \begin{bmatrix} 0 & -a_z & a_y \\ a_z & 0 & -a_x \\ -a_y & a_x & 0 \end{bmatrix} = \begin{bmatrix} a_z^2 + a_y^2 & -a_x a_y & -a_x a_z \\ -a_x a_y & a_x^2 + a_z^2 & -a_y a_z \\ -a_x a_z & -a_y a_z & a_x^2 + a_y^2 \end{bmatrix} = \boldsymbol{a}^T\boldsymbol{a}.\boldsymbol{I} - \boldsymbol{a} \otimes \boldsymbol{a}^T \tag{14}$$

where $\otimes$ means a *tensor product* (result is a matrix).

Now, we can write using:

$$\boldsymbol{M}_p = \boldsymbol{Q}_p^T\boldsymbol{Q}_p \qquad \boldsymbol{R} = \sum_{p=1}^{2N} \boldsymbol{M_p} \tag{15}$$

The matrix $\boldsymbol{M}_p$ can be expressed as:

$$\boldsymbol{M}_p = \boldsymbol{x}_p^T\boldsymbol{x}_p \ .\boldsymbol{I} - \boldsymbol{x}_p \ \otimes \boldsymbol{x}_p^T \tag{16}$$

Then, the matrix $\boldsymbol{R}$ can be expressed as:

$$\boldsymbol{R} = \sum_{p=1}^{N} \left(\boldsymbol{x}_p^T\boldsymbol{x}_p \ .\boldsymbol{I} - \boldsymbol{x}_p \ \otimes \boldsymbol{x}_p^T\right) = \sum_{p=1}^{N} \boldsymbol{x}_p^T\boldsymbol{x}_p \ .\boldsymbol{I} - \sum_{p=1}^{N} \boldsymbol{x}_p \ \otimes \boldsymbol{x}_p^T = \xi\boldsymbol{I} - \boldsymbol{\Omega} \tag{17}$$

where $N$ is a number of points and

$$\xi = \sum_{p=1}^{N} \boldsymbol{x}_p^T\boldsymbol{x}_p \qquad \boldsymbol{\Omega} = \sum_{p=1}^{N} \boldsymbol{x}_p \ \otimes \boldsymbol{x}_p^T \tag{18}$$

An extreme of the distance $D$ is given as:

$$\frac{\partial D}{\partial \boldsymbol{s}} = \boldsymbol{0} \tag{19}$$

Now

$$\begin{aligned} \frac{\partial D}{\partial \boldsymbol{s}} &= \frac{\partial}{\partial \boldsymbol{s}} \sum_{p=1}^{N} \frac{\left(\boldsymbol{x_p} \times \boldsymbol{s}\right)^T\left(\boldsymbol{x_p} \times \boldsymbol{s}\right)}{\boldsymbol{s^T s}} = \sum_{p=1}^{N} \frac{\partial}{\partial \boldsymbol{s}} \frac{\left(\boldsymbol{x_p} \times \boldsymbol{s}\right)^T\left(\boldsymbol{x_p} \times \boldsymbol{s}\right)}{\boldsymbol{s^T s}} \\ &= \sum_{p=1}^{N} \frac{\partial}{\partial \boldsymbol{s}} \frac{\left(\boldsymbol{Q}_p\boldsymbol{s}\right)^T\left(\boldsymbol{Q}_p\boldsymbol{s}\right)}{\boldsymbol{s^T s}} = \sum_{p=1}^{N} \frac{\partial}{\partial \boldsymbol{s}} \frac{\boldsymbol{s}^T\boldsymbol{Q}_p^T\boldsymbol{Q}_p\boldsymbol{s}}{\boldsymbol{s^T s}} = \sum_{p=1}^{N} \frac{\partial}{\partial \boldsymbol{s}} \frac{\boldsymbol{s}^T\boldsymbol{M}_p\boldsymbol{s}}{\boldsymbol{s^T s}} \\ &= \frac{\partial}{\partial \boldsymbol{s}} \frac{\boldsymbol{s^T R\ s}}{\boldsymbol{s^T s}} = \boldsymbol{0} \end{aligned} \tag{20}$$

Now, for the extreme we get:

$$\frac{\partial D}{\partial \boldsymbol{s}} = \frac{2\boldsymbol{Rs}\,.\boldsymbol{s}^T\boldsymbol{s} - \boldsymbol{s}^T\boldsymbol{R}\,\boldsymbol{s}\,.\ 2\boldsymbol{s}}{(\boldsymbol{s}^T\boldsymbol{s})^2} = \boldsymbol{0} \tag{21}$$

We have to solve a system of equations; note it is a ***vector*** form:

$$(\boldsymbol{s}^T\boldsymbol{s})\,.\boldsymbol{R}\,\boldsymbol{s} - (\boldsymbol{s}^T\boldsymbol{R}\,\boldsymbol{s})\,.\ \boldsymbol{s} = \boldsymbol{0} \qquad \text{i.e.} \qquad (\boldsymbol{s}^T\boldsymbol{s}\,.\boldsymbol{R} - \boldsymbol{s}^T\boldsymbol{R}\,\boldsymbol{s}\,.\boldsymbol{I})\,.\ \boldsymbol{s} = \boldsymbol{0} \tag{22}$$

where $\boldsymbol{I}$ is an identical matrix.

The system of equations Eq.(22) can be simplified as follows:

$$\begin{aligned} (\boldsymbol{s}^T\boldsymbol{s})\,.\boldsymbol{R}\,\boldsymbol{s} - (\boldsymbol{s}^T\boldsymbol{R}\,\boldsymbol{s})\,.\ \boldsymbol{s} &= (\boldsymbol{s}^T\boldsymbol{s})\,.(\xi\boldsymbol{I} - \boldsymbol{\Omega})\,.\boldsymbol{s} - (\boldsymbol{s}^T(\xi\boldsymbol{I} - \boldsymbol{\Omega})\,\boldsymbol{s})\,.\boldsymbol{s} \\ &= \xi(\boldsymbol{s}^T\boldsymbol{s})\,.\boldsymbol{s} - (\boldsymbol{s}^T\boldsymbol{s})\,.\boldsymbol{\Omega}\,\boldsymbol{s} - (\xi\boldsymbol{s}^T\boldsymbol{s} - \boldsymbol{s}^T\boldsymbol{\Omega}\,\boldsymbol{s})\,.\boldsymbol{s} \\ &= -(\boldsymbol{s}^T\boldsymbol{s})\,.\ \boldsymbol{\Omega}\,\boldsymbol{s} + \boldsymbol{s}^T\boldsymbol{\Omega}\,\boldsymbol{s}\,.\boldsymbol{s} = \boldsymbol{0} \end{aligned} \tag{23}$$

i.e.

$$(\boldsymbol{s}^T\boldsymbol{s})\,.\ \boldsymbol{\Omega}\,\boldsymbol{s} - \boldsymbol{s}^T\boldsymbol{\Omega}\,\boldsymbol{s}\,.\boldsymbol{s} = \boldsymbol{0} \tag{24}$$

Due to the implicit formulation, we can divide the Eq.(24) by any constant $\xi \neq 0$. Let us set $\xi = \|\boldsymbol{s}\|$ and divide the Eq.(24) by $\xi^2$, i.e. we will use a normalized vector $\|\boldsymbol{s}\|$=1, i.e. $\boldsymbol{s}^T\boldsymbol{s} = 1$.

Then we can write:

$$\boldsymbol{\Omega}\, \boldsymbol{s} - \boldsymbol{s}^T \boldsymbol{\Omega}\, \boldsymbol{s} \,.\, \boldsymbol{s} = \boldsymbol{0} \tag{25}$$

i.e.

$$(\boldsymbol{\Omega} - (\boldsymbol{s}^T \boldsymbol{\Omega}\, \boldsymbol{s})\, \boldsymbol{I}\,) \,.\, \boldsymbol{s} = \boldsymbol{0} \tag{26}$$

where: $(\boldsymbol{s}^T \boldsymbol{s})$, $\boldsymbol{s}^T \boldsymbol{\Omega}\, \boldsymbol{s}$ are scalar values, $\boldsymbol{\Omega}$ is a matrix, $\boldsymbol{s}$ is a vector.

The "unspecified" value $\xi \neq 0$ has no influence as the vector $\boldsymbol{s}$ is a directional vector and has no influence to the computed distance $D$ value. We obtain the required line $\pi$ by solving the system Eq.(26), as the $\boldsymbol{x}_A = \boldsymbol{0}$. Note, that if the given points were shifted, the point $\boldsymbol{x}_A$ is to be shifted back as well.

The Eq.(26) is not actually dependent on a given data dimensionality, but it is dependent on a dimensionality of a geometric primitive we are mapping to.

## CONCLUSION

This paper presents a new formulation for the Total (Orthogonal) Least Square Error method for the case when points are given in the $E^3$ space and a line $\pi$, fitting the TLSE criteria, is to be found. The derivation of the final formula is simple. However there is a question how to solve it efficiently. A typical application of the presented approach is when some samples occur in 3 dimensional parametric space and theoretically parameters values should be linearly depended, i.e. all samples should lie on a line in 3 dimensional space.

Generalization of the presented approach to the $d$-dimensional case is straightforward. On the other hand, cases when data are given in the $n$-dimensional space and hyperplane in $k$-dimensional space is searched is to be solved as well. In future work the given approach is to be extended to general case, when data are given in the $n$-dimensional space and hyperplane in $k$-dimensional space, where $k < n$, is searched using geometry algebra approach, e.g. how to compute the TLSE in 4-dimensional case if mapping is to be made to a plane in the 4-dimensional case.

## ACKNOWLEDGMENT

The author thanks to students and colleagues at the University of West Bohemia for recommendations, constructive discussions, and hints that helped to finish the work. Many thanks belong to the anonymous reviewers for their valuable comments and suggestions that improved this paper significantly. This research was supported by the Ministry of Education of the Czech Republic – project No.LH12181.